\documentclass[10pt]{amsart}
\usepackage[english]{babel}
\usepackage{amsmath} \usepackage{amsthm}
\usepackage{amssymb}
\usepackage{amscd}
\usepackage{tabularx}
\newtheorem{lem}{Lemma}[section]
\newtheorem{prop}[lem]{Proposition}
		\newtheorem{thm}[lem]{Theorem}
		
		\newtheorem{cor}[lem]{Corollary}
\theoremstyle{definition}
\newtheorem{Def}[lem]{Definition}
\newtheorem{ese}[lem]{Example}
\theoremstyle{remark}
		\newtheorem{rem}{Remark}
		\DeclareMathOperator{\rank}{rank}
		\DeclareMathOperator{\Hom}{Hom}

\def\sist{\left\{\begin{array}{l}}
\def\sistt{\end{array}\right.}
\def\bb{\mathbb}
\def\la{\langle}
\def\ra{\rangle}
\def\phi{\varphi}

\begin{document}
\title{On projective varieties of
dimension $n+k$ covered by $k$-spaces}
\author[E.\ Mezzetti, O.\ Tommasi] {E.\ Mezzetti $^1$, O.\ Tommasi}
  \address{Dipartimento di Scienze
Matematiche, Universit\`a di Trieste, Via Valerio 12/1, 34127 Trieste, Italy}
\email{mezzette@univ.trieste.it}
\address{Department of Mathematics, University of Nijmegen, Toernooiveld, 
6525 ED Nijmegen,  The Netherlands}
\email{tommasi@sci.kun.nl}
\date{}
\subjclass{ Primary 14N20, 14J30; Secondary 53A20, 14D05, 14M05}
\keywords{Focal scheme, second fundamental form, Gauss map, tangent
space to Grassmannian, 
ruled threefolds}
\thanks{$^1$ Supported by University of Trieste (fondi 60\%), MURST,
project \lq\lq Geometria sulle variet\`a algebriche'' and by  INdAM, 
project \lq\lq Cohomological
properties of projective schemes''}
\headsep=40pt
\footskip=45pt
\headheight 20pt

\begin{abstract}  We study families of linear spaces in 
projective space 
whose union is a proper subvariety $X$ of the expected dimension.
We establish relations between configurations of focal points and existence
or non-existence of a fixed tangent space to $X$ along a
general element of the family. We apply our results to the 
classification of ruled $3$-dimensional varieties. \end{abstract}
\maketitle

\section*{Introduction}
Since the publication of \cite{GH79} there has been a renewal of
interest in the study of differential geometric properties of 
algebraic varieties.
The bases of this study are to be found in classical
works, such as several papers by C. Segre (particularly  
\cite{S:sup} and \cite{S:fuochi}). 
There, topics such as the second fundamental
form of projective varieties, varieties with degenerate Gauss mapping and
in general varieties ruled by linear subspaces are introduced and
discussed. 
Recently, contributions on these topics have been given by
Akivis, Goldberg, Landsberg, Rogora
(\cite{AG},\cite{L99},\cite{AGL},\cite{R}). These papers highlight
the importance of the study of the {\em focal scheme}.

The foci are
a classical tool for families of linear spaces (see \cite{S:fuochi}).
In modern algebraic geometry it has been reformulated by means of the
focal diagram in the paper of Ci\-li\-ber\-to and Ser\-ne\-si (\cite{CS})
and has been applied  to the study of congruences of lines
(\cite{ABT},\cite{Arr},\cite{DePoi}).

In this paper we will
deal with families of linear spaces that generate proper subvarieties of the
expected dimension in the
projective space. For instance, let us 
consider a family
$B$ of $k$-spaces in the projective space $\bb{P}^N$,
the variety $X$ ruled by $B$, and assume $\dim B=n$, $\dim X= n+k<N$. 
Then, we  will take into consideration  
the relationship between the existence and
the properties of the focal scheme on a general space of $B$,
and the existence of spaces of dimension $\leq n+k$ tangent to $X$ along
a general space of $B$. A complete description of this
relationship for a family of lines will be given in the following theorem:

\begin{thm}\label{thm:C}
Let $B\subset\mathbb{G}(1,N)$ be a family of lines in $\mathbb{P}^N$
of dimension $n,\ n\leq N-2$. Suppose that the union of
the lines belonging to $B$ is an algebraic variety $X$ of dimension $n+1$.
Then, for all $k$ in the range $0\leq k\leq n$, the following are 
equivalent:\\
(i) the focal locus on the general element $r\in B$
has length $k$;\\
(ii) X has a fixed tangent  $\mathbb{P}^{k+1}$ along every general
$r\in B$. \end{thm}

There is an analogue 
to Theorem \ref{thm:C} for varieties with degenerate
Gauss mapping:

\begin{thm}\label{thm:B}
Let $B$ be a family of linear subspaces of $\mathbb{P}^N$ of dimension
$k$, and denote by $n$ the dimension of $B$. Suppose that the union of
the $k$-planes of the family $B$ is an algebraic variety
$X\subset\mathbb{P}^N$ of dimension $n+k<N$. Then the following are
equivalent:\\
(i) the tangent space to $X$ is constant along general
elements of $B$;\\
(ii) for all $\Lambda$ belonging to an open set
of $B$ the focal subvariety of $B$ is a hypersurface of $\Lambda$ of
degree $n$; otherwise all points of $\Lambda$ are focal.
\end{thm}

We will apply Theorem \ref{thm:C} and Theorem \ref{thm:B} to the
study of ruled varieties of dimension 3.
Our results  comprise and
complete what is shown in previous papers, such as
\cite{GH79}, \cite{R}, \cite{AGL}. It should be noted, however, that
the result in
\cite{GH79} about varieties with degenerate Gauss mapping
is not precisely stated, and that
\cite{R} considers only necessary conditions and not sufficient ones.
We will give the classification of
threefolds with a tangent 2-plane constant along lines in Theorem \ref{thm:c1},
and that of threefolds with
degenerate Gauss mapping  in Theorem \ref{thm:c2}.

\begin{thm}\label{thm:c1}
Let $B$ be a surface in the Grassmannian $\mathbb{G}(1,N)$, 
with $N\geq4$. 
Suppose that the union of the lines belonging to $B$ is an algebraic 
variety $X$
of dimension 3, and that the Gauss image of $X$ has
dimension 3. Then, along a general line of $B$ there is a fixed tangent
  2-plane not contained in $X$ if and only if $X$ is one of
the following:\begin{enumerate}
\item a union of lines, all tangent to a surface $S\subset\bb{P}^N$,
  whose direction  at the tangency point 
is not in general a conjugate direction for the second fundamental
form of $S$;
\item the union of a one-dimensional family of 
2-dimensional cones,
whose vertices sweep a curve.
\end{enumerate}\end{thm}
\begin{thm}\label{thm:c2}
Let $X$ be a variety of dimension 3 with degenerate Gauss mapping.
Then, one of the following holds 
\begin{enumerate}
\item\label{G2} the Gauss image of $X$ has dimension 2, and $X$ is one of
the following:
\begin{enumerate}
\item a union of lines bitangent to a surface;
\item there are two surfaces such that $X$ is a union of lines tangent to
  both;
\item a union of lines tangent to a surface, and meeting a fixed curve;
\item the union of asymptotic tangent lines of a surface;
\item the join of two curves;
\item the variety of secant lines of a curve;
\item a band (see Definition \ref{band});
\item the cone over a surface, with a point as vertex.
\end{enumerate}
\item\label{G1} the Gauss image of $X$ has dimension 1,
and $X$ is  built up by a composite construction of cones
and varieties ruled by osculating spaces over some curve. \end{enumerate}
All these cases are possible, and each of them always
represents a class of varieties with degenerate Gauss mapping.
\end{thm}
The plan of the paper is as follows: In Section
\ref{sec:foci}
we  introduce the notion of foci for a family of linear
spaces and
we give the interpretation of foci in terms of tangent
spaces to the Grassmannian.

In Section \ref{dim3} we prove the two classification theorems for 
ruled varieties of dimension 3. We prove moreover that all 
surfaces $S$ appearing in
cases (a)-(d) of Theorem  \ref{thm:c1} are not general, but must satisfy 
the condition that the osculating space to $S$ at a general point 
has dimension at most 4. 

In Section \ref{s:rette}
we  prove Theorem
\ref{thm:C}  and establish
the properties of the focal locus in the case of varieties
ruled by lines. In the last section we  consider varieties
ruled by subspaces of higher dimension.
We  prove by means of an example
that Theorem \ref{thm:C} cannot be extended to a family
of subspaces of dimension $\geq2$. 
However, Theorem \ref{thm:B} shows that
a description is still possible for
varieties with degenerate Gauss mapping. This result has 
already
been proved by Akivis and Goldberg in \cite{AG} with differential
 geometry techniques; we  give now an algebraic proof of it.

\section*{Notation}
We will study projective algebraic varieties over the complex field
or, more generally,
over an
algebraically closed field $\mathbb{K}$ with 
$char\ \mathbb{K} =0$.\\
$V$ will denote a linear space of dimension $N+1$ over $\mathbb{K}$,
and $\mathbb{P}^N=\mathbb{P}(V)$ the projectivization of $V$. Analogously,
$\bb{A}^{N+1}=\bb{A}(V)$
will denote the affine space associated to
$V$.\\ If $\Lambda\subset \bb{P}^N$ is a projective linear subspace,
$\hat{\Lambda}\subset V$ will denote the linear subspace associated
to $\Lambda$ such that $\Lambda=\bb{P}(\hat{\Lambda})$.\\
$[v] \in \mathbb{P}(V)$ will denote the point of $\bb{P}^N$ corresponding
to the equivalence class of $v\in V\smallsetminus\{0\}$.\\
$T_xX$ will denote the Zariski tangent space to the variety $X$ at
  its point $x$,
while we will denote by $\mathbb{T}_xX\subset\bb{P}^N$ the embedded
tangent space to $X$ at $x$.\\
$G(h,V)$ will denote the Grassmannian variety of linear subspaces of
dimension $h$ in $V$. $\mathbb{G}(k,N)$ will denote the Grassmannian 
of projective
subspaces of dimension $k$ of $\bb{P}^N$. We will use
the same symbol to denote the points of the Grassmannian and the
corresponding linear subspaces.

\section{Focal diagram}\label{sec:foci}
Let $B\subset\bb{G}(k,N)$ be a family of dimension $n$ of $k$-spaces 
in $\bb{P}^N$.
Denote by $B'$ a desingularization of $B$ and by $\mathcal{I}$ 
the incidence correspondence of $B'$, with the natural projections
$$\begin{array}{ccccc} \;B'&\xleftarrow[\mspace{50mu}]{p_1}&B'
\times\bb{P}^N&\xrightarrow[\mspace{50mu}]{p_2}&\bb{P}^N\\[-6pt] 
&&\cup &&\\[-6pt]
B&\xleftarrow[\mspace{50mu}]{g}&\mathcal{I}&\xrightarrow[\mspace{50mu}]{f}&\;
\mathbb{P}^N.\end{array}$$
In what follows, we will restrict ourselves to families $B$ such that 
the image of $f$ (i.e. the union of the lines belonging to $B$)
is a variety $X$ of dimension $n+k$. This is the same as assuming
the general fibre of $f: \mathcal{I}\longrightarrow X$ to be finite.
\begin{Def}
A point $x\in X$ is said to be a {\em fundamental point} of the family
$B$ if the fibre $f^{-1}(x)$ has positive dimension.\\
This condition defines a closed subset of $X$ called the {\em 
fundamental locus} $\Phi$ of $B$.
\end{Def}
On the basis of this set-up, we can construct a
commutative diagram of exact sequences, called the {\em focal diagram} of $B$:
$$\begin{CD}
@. @. 0 @.\\
@. @. @VVV\\
@. @. \left( p_1^*(\mathcal{T}_{B'})\right)\,_{|\mathcal{I}}@>\chi>>
\mathcal{N}_{\mathcal{I}|B' \times \mathbb{P}^N}\\
@. @. @VVV @|\\
0 @>>> \mathcal{T}_{\mathcal{I}} @>>> \mathcal{T}_{B' \times
\mathbb{P}^N}{}_{|\mathcal{I}} @>>> \mathcal{N}_{\mathcal{I}|B' 
\times \mathbb{P}^N} @>>> 0\\
@. @VdfVV @VVV\\
@. f^*(\mathcal{T}_{\mathbb{P}^N}) @=
\left( p_2^*(\mathcal{T}_{\mathbb{P}^N})\right)\,_{|\mathcal{I}}@.\\ 
@. @. @VVV\\
@. @. \;0.
\end{CD}$$
The focal diagram is built up by crossing
the exact sequence
defining the normal sheaf to $\mathcal{I}$ inside $B'\times\bb{P}^N$
with the sequence (restricted to $\mathcal{I}$) expressing the
tangent sheaf of the product variety $B'\times\bb{P}^N$ as a product
of tangent sheaves.

\begin{Def}
The map denoted by $\chi$ in the focal diagram is
called the {\em characteristic map} of the family $B$.\\
For every $\Lambda\in B_{ns}$ the restriction of $\chi$ to $g^{-1}(\Lambda)$
is called the {\em characteristic map of $B$ relative to $\Lambda$};
it lies in the following diagram:
$$\begin{CD}
  \chi(\Lambda):\ @. T_\Lambda B'
\otimes \mathcal{O}_\Lambda @>>> \mathcal{N}_{\Lambda|\mathbb{P}^N}\\
\scriptstyle @. \wr@| \wr@|\\
\textstyle @. \mathcal{O}^{m}_\Lambda @>>> \mathcal{O}^{N-k}_\Lambda (1). \end{CD}$$
\end{Def}
\begin{Def}
The condition
$$\rank{\chi} (\Lambda,x) < \min \{\rank(\left( p_1^*(\mathcal{T}_{B'})
\right)\,_{|\mathcal{I}}), \rank(\mathcal{N}_{\mathcal{I}|B' \times 
\mathbb{P}^N})\}$$
defines a closed subscheme $V(\chi)\subset\mathcal{I}$ which will
be called the \emph{subscheme of first order foci} (or, simply, the
{\em focal subscheme}) of the family $B$. Analogously, $F = f
(V(\chi))$ is called the \emph{locus of first order foci}, or the
{\em focal locus} of $B$ in $\mathbb{P}^N$.
\end{Def}
By the commutative property of the focal diagram, the focal
locus has a double interpretation. Indeed, the kernel of $\chi$ 
and the kernel of $df$ must coincide (as subsheaves of 
$\mathcal{T}_{B
\times \mathbb{P}^N}{}_{|\mathcal{I}})$). Then the focal locus
is the ramification locus of $f$. As a consequence, the fundamental
locus is contained in the focal locus. These considerations can be rephrased
by the following proposition.
\begin{prop}\label{equiv}
The following
  are equivalent:\\
1. the rank of $\chi$ is maximal;\\
2. the rank of $df$ is maximal;\\
3. $V(\chi)$ is a closed proper subscheme of $\mathcal{I}$. \end{prop}

We have assumed that the union of the $k$-spaces belonging to $B$
is a variety $X$ of dimension $n+k$. By Proposition \ref{equiv},
this implies that a general
point on a general space of $B$ is not a focus. Nevertheless,
  some particular spaces of $B$ can be contained in the focal locus:
they are called {\em focal spaces}. \vskip 6pt
The characteristic map is
closely connected with the structure of the tangent space to
the Grassmannian variety as a space of homomorphisms (see \cite{H92}).
Let $B$ be a subvariety of $\bb{G}(k,N)$. We can identify this Grassmannian
with the Grassmannian of linear subspaces of dimension $k+1$ of
$V$, $G(k+1,V)$. Then, by associating to each $\Lambda\in B$ the affine
  cone $\hat{\Lambda}\subset\bb{A}(V)= \bb{A}^{N+1}$,
  we can construct a new incidence correspondence 
$\mathcal{I}'\subset B'
\times\bb{A}^{N+1}$ and projections $f',g'$:
$$\begin{array}{ccccc} \;B'&\xleftarrow[\mspace{50mu}]{q_1}&B'
\times\bb{A}^{N+1}
&\xrightarrow[\mspace{50mu}]{q_2}&\bb{A}^{N+1}\\[-6pt] &&
\mspace{-18mu}\cup&&\\[-6pt]
B&\xleftarrow[\mspace{50mu}]{g'}&\mspace{-18mu}\mathcal{I'}\ &
\xrightarrow[\mspace{50mu}]{f'}&\;
\mathbb{A}^{N+1}.\end{array}$$
Considering
$B$ as a family of subspaces in $\bb{A}^{N+1}$ yields an affine 
version of the focal diagram:
$$\begin{CD}
@. @. 0 @.\\
@. @. @VVV\\
@. @. \left( q_1^*(\mathcal{T}_{B'})\right)\,_{|\mathcal{I}'}@>\chi'>>
\mathcal{N}_{\mathcal{I}'|B' \times \mathbb{A}^{N+1}}\\
@. @. @VVV @|\\
0 @>>> \mathcal{T}_{\mathcal{I}'} @>>> \mathcal{T}_{B' \times
\mathbb{A}^{N+1}}{}_{|\mathcal{I}'} @>>> \mathcal{N}_{\mathcal{I}'|B' \times
\mathbb{A}^{N+1}} @>>> 0\\
@. @Vdf'VV @VVV\\
@. {f'}^*(\mathcal{T}_{\mathbb{A}^{N+1}}) @=
\left( q_2^*(\mathcal{T}_{ \mathbb{A}^{N+1}})\right)\,_{|\mathcal{I}'}@.\\
@. @. @VVV\\
@. @. \;0.
\end{CD}$$
As in the projective case, we can define the characteristic
map $\chi'$ relative to $\Lambda$, a non-singular element of $B$,
$$\chi'(\Lambda):\ T_\Lambda B\otimes\mathcal{O}_{\hat{\Lambda}}
\longrightarrow \mathcal{N}_{\hat{\Lambda}|\bb{A}^{N+1}}.$$
If we compare the definition of the focal diagram and
the characterization of $T_\Lambda B$ as a space of homomorphisms,
we easily get the following proposition.
\begin{prop}\label{prop:h&f}
Let $\Lambda$ be a non-singular point of $B\subset\bb{G}(k,N)=G(k+1,V)$.
Let us consider $T_\Lambda B$ as a linear subspace of $T_\Lambda\bb{G}(k,N)
\cong\Hom(\hat{\Lambda},V/\hat{\Lambda})$. Then the
characteristic map $\chi'$ relative to $\Lambda$, considered
as a morphism of vector bundles, for all $v\in \hat{\Lambda}$ associates to
$\eta\in T_\Lambda B$ the normal vector $\eta (v)$.
\end{prop}
The projectivization of the usual characteristic map $\chi$
coincides with that of the affine version $\chi'$, so we have:
\begin{cor}\label{cor:h&f}
Let $\Lambda$ be a non-singular point of $B\subset\bb{G}(k,N)$.
Let us consider $T_\Lambda B$ as a linear subspace of
$T_\Lambda\bb{G}(k,N)\cong\Hom(\hat{\Lambda},V/\hat{\Lambda})$.
Then the projectivization of the characteristic map $\chi$ relative
to $\Lambda$, considered as a morphism of vector bundles, for all
$p\in\Lambda$ associates to $[\eta]\in \bb{P}(T_\Lambda B)$
the point $[\eta (v)]\in\bb{P}(V/\hat{\Lambda})$, where $v\in V$ is such that
$[v]=p$. \end{cor}

This corollary yields  an interpretation of
focal points, which is particularly clear in the case of a family of lines.

\begin{rem}\label{fuochimult}
Consider a variety $B\subset\bb{G}(1,N)$ and a general line
$r\in B$. Then, the foci on $r$ are the points $p=[v]$ such
that $v\in\ker\eta$
for a non-trivial $\eta\in T_rB$. Since under our hypotheses 
the rank of the general $\eta\in T_rB$ is 2, the existence 
of focal points depends on the existence of rank 1 
homomorphisms in $T_rB$. 
  Focal points with multiplicity represent a special case. A point
$p=[v]\in r$ is a focal point of multiplicity $\geq2$ if and only
if there exist two linearly independent tangent vectors
$\eta_1,\eta_2\in T_rB$ verifying
$$\begin{array}{ll} \eta_1(v)=0,&Im\ (\eta_1)\neq0,\\
\eta_2(v)\in Im\ (\eta_1),& Im\ (\eta_2)\neq Im\ (\eta_1),\end{array}$$
since the condition on multiplicity is that the composition of the 
characteristic
map relative to $r$ with the natural map $V/\hat{r}
\longrightarrow (V/\hat{r})/Im\ (\eta_1)$ has not maximal rank.
Iteration of this construction provides the characterization
for focal points of higher multiplicity.  
\end{rem}

By means of the Pl\"ucker embedding, we can consider
the embedded tangent
  space to
$B\subset\bb{G}(1,N)$ at a point $\Lambda$. In the case of lines,
there is a connection
between the existence of focal points on $r\in B$
and the existence of a line in $\bb{T}_rB\cap\bb{G}(1,N)$.
\begin{prop}
Le $B\subset\bb{G}(1,N)$ be a family of lines. Let $r$
be a general element of $B$.
Suppose $r$ is not focal: then, there is a bijection
between the focal points on $r$ and
the lines in the intersection of the Grassmannian
$\bb{G}(1,N)$ with $\bb{T}_rB$ (embedded in $\bb{P}(\bigwedge^2V)$). \end{prop}
\proof
We know that $[v]\in r$ is a focal point if and only if $v\in\ker\eta$,
where $\eta\in T_rB$ has rank 1. With a simple computation,
it is possible to prove that if a homomorphism $\eta\in T_r\bb{G}(1,N)$ has
rank 1 then the pencil of lines passing through $\bb{P}(\ker\eta)$ and
lying in $Im\ (\eta) \oplus r$ is a line contained in the
intersection of $\bb{T}_rB$ with the Grassmannian.
The converse is also true:
if there is a line in the intersection, then we can find a
homomorphism $\eta$ of rank 1 and hence a focal point.
\qed

\section{Varieties of dimension 3}\label{dim3}
We will apply the study of the focal locus to the specific
problem
of classifying ruled varieties of dimension 3 with
degenerate tangential properties. More precisely, we will consider:
\begin{enumerate}
\item\label{a} varieties ruled by lines with a
constant tangent 2-plane along
any line of the ruling;
\item varieties ruled by lines with a constant tangent space of dimension 3
along every line;
\item varieties ruled by planes with a constant tangent space of dimension
3 along every plane.
\end{enumerate}
Note that the last two cases yield the classification
of varieties of
dimension 3 with degenerate Gauss mapping. In the proofs we will
use also some results to be proved in Sections 3 and 4.

The classical references for our approach to classification 
are the works of C. Segre. Particularly, a classical proof of
the classification of case \ref{a} can be found in \cite{S:fuochi}.
The classification of varieties of dimension 3 with degenerate
Gauss mapping has
already been presented recently in \cite{R} and \cite{AGL}.
In both papers, the classification is based on the study of
the focal scheme
of the family of fibres of the Gauss map, but there
is no distinction between
strict focal locus and (total) focal locus
(see Definition \ref{d:strict}). In \cite{R} the
classification is outlined 
without 
a study of the second fundamental form of focal surfaces.
Therefore, there is no description of how to construct a
variety with degenerate Gauss mapping. In \cite{AGL}
one of the cases (that of bands) is not
completely solved.
\vskip 12 pt
In what follows, the concept of {\em conjugate directions}
for the second fundamental form will naturally arise. We will denote
by $$II_y: T_yY\otimes T_yY\longrightarrow N_yY$$
the {\em second fundamental form} of a variety $Y$ at a non-singular point
$y$ (for the definition, see \cite{H92}). It is a symmetric bilinear 
form, so it can be
interpreted as a linear system of quadrics
$|II_y|$ in
$\bb{P}(T_yY)$. The dimension of the linear system is linked with the
dimension of the second osculating space to $Y$ in $y$,
$T^{(2)}_yY$, by the relation
$$\dim|II_y|=\dim T^{(2)}_yY-\dim T_yY-1.$$
\begin{Def}
Let $Y\subset\bb{P}^N$ be a variety, and $y$ be a non-singular point of it.
Then two tangent vectors $v,w\in T_yY$ are said to represent {\em 
conjugate directions}
at $y$ if $II_y(v,w)=0$. This means that the
points $[v],[w]$ are conjugate with respect to all quadrics in 
$|II_y|$.\\ If there is a
selfconjugate tangent vector, its direction is
called an {\em asymptotic direction} at $y$.
\end{Def}
The existence of conjugate directions at every non-singular point is
not a general fact. It is related to the dimension of the second
osculating space to the variety at the general point. We are
interested in the study of conjugate directions for surfaces.
For general surfaces at general points the dimension of the second
osculating space is 5. In this case, at a general point there
are no conjugate directions. Conjugate directions exist only if the
dimension of the second osculating space is $\leq4$.
If the dimension is 3, every direction possesses a conjugate 
direction. It is well known
that in $\bb{P}^N$, $N\geq4$, this property holds only for 
developable surfaces, i.e. cones and varieties
swept out by tangent lines to a curve.
\begin{Def}[\cite{S:sup}]
A surface is called a $\varPhi$ surface if and only
if the dimension of its second osculating space at the general point is 4.
\end{Def}

\begin{prop}
For a surface $S\subset\bb{P}^N$, $N\geq5$, the following properties are
equivalent:
\\(i) $S$ is a $\varPhi$ surface; 
\\(ii) at a general point of $S$ there is exactly one couple of conjugate
directions, or one asymptotic direction;
\\(iii) the union of the tangent planes to $S$ is a
variety of dimension $4$
with tangent space fixed along those planes.
\end{prop}
\proof The equivalence of the first two properties is a
consequence of the fact that a linear system of quadrics in $\bb{P}^1$
admits exactly one couple of conjugate points if and only if its 
dimension is 1.
Let now $S\subset\bb{P}^N$ be a surface: let us denote by $V$ the
closure in $\bb{P}^N$ of the union of the tangent planes to $S$ at
its non-singular points. Then the dimension of $V$ is 4 if and only if
$S$ is not a developable surface or a plane. Let us consider the
osculating space $\bb{T}^2_pS$ to $S$ at a general point $p$, embedded
in $\bb{P}^N$. Using a local parametric representation of $S$, it is
  easy to show the following equality:
$$\bb{T}^2_pS=\overline{\bigcup_{q\in\bb{T}_pS\cap V_{ns}}\bb{T}_qV}.$$
This implies that the tangent space to $V$ is constant along planes
if and only if the dimension of the osculating space to $S$ equals the
dimension of $V$. 
Hence, the equivalence of (i) and (iii) is established. \qed

\begin{rem}
The general situation for the union of tangent planes
to a surface is that the fibres of the Gauss map are 1-dimensional.
\end{rem}

When we have a $\varPhi$ surface $S$, we can always construct an 
irreducible family $\Sigma$ of dimension 2, whose elements are 
lines tangent to $S$, such that for each general point $p\in S$
there is exactly one line of $\Sigma$ tangent to $S$ at $p$, and
moreover its tangent direction at $p$ is conjugate to some (other) 
tangent direction.
In this
case, we will say that the lines of $\Sigma$ {\em admit a conjugate 
direction} on $S$. If at
  the general point of $S$ the
two conjugate directions coincide, i.e. there is an asymptotic
direction, then the lines of $\Sigma$ are called
{\em asymptotic lines} on $S$.

\vskip 12pt
Let us consider case \ref{a} first. In this case,
we have a 3-dimensional variety $X$ which is covered by the lines
belonging to a surface $B$ in the
Grassmannian $\bb{G}(1,N)$, such that $X$ has a constant tangent
plane along a general line of $B$.
  A classification of these families is provided by
  Theorem \ref{thm:c1} (see also \cite{M}).

\vskip 6pt
\par{\em Proof of Theorem \ref{thm:c1}. }
We can apply Theorem \ref{thm:C} to the family $B$.
The existence of the tangent plane implies
then that on the general line belonging to $B$
there exists one focal point
(with multiplicity 1).
This focus cannot be a
fixed point $p$.
In this case there would be
a $2$-dimensional subfamily of lines of $B$ passing through $p$, and $p$
would be a focal point of multiplicity 2, which is not allowed.
Then, considering the closure of the union of the focal points on
such lines (the {\em strict} focal locus, in the terminology to be
introduced in \S 3), we get two possibilities: we can obtain a surface
$S$, or a curve $C$. In the former case, the first part of the claim
follows from Theorem \ref{thm:A}. The exception considered in our
statement is necessary in order to exclude varieties with degenerate
Gauss mapping, as we will see later. In the latter case, $C$ lies
in the fundamental locus of $B$, which yields the second part of the
  claim.
By a direct calculation, we can 
check that the union of tangent
lines to a surface $S$ has a constant tangent plane along a
general line $r$. This plane is the tangent plane to $S$ at the
  point of tangency of $r$. Analogously, the fixed tangent plane along the
lines of a cone is contained in the tangent space to the union of cones. \qed

\begin{rem}
In the hypotheses of Theorem \ref{thm:c1} we have excluded the (trivial)
case of varieties $X$ ruled by planes. In this case, the Fano variety of
lines has dimension $>2$, but it is always
possible to find a subvariety $B$ of it, with $\dim B=2$, such that
the lines of $B$ cover $X$. There are two ways of
constructing $B$. We can choose a unisecant curve $C$ to the family of
planes and consider for every plane the pencil of lines
with center the corresponding point of $C$. The points of $C$ are fundamental
points of $B$, and in general they are not singular points of $X$.
Note that a general ruled surface in $\bb{G}(1,N)$ gives an example of
this situation. Otherwise, inside every plane we can fix a  curve
(varying algebraically with the plane) and consider the family of its
tangent lines. Also in this case the focal points are not in general
singular for $X$. In fact, in both cases the focal points have no real
geometric meaning for $X$. \end{rem}

We will prove now Theorem \ref{thm:c2}, giving the classification of
varieties of dimension 3 with degenerate Gauss
mapping.

\vskip 6pt
\par{\em Proof of Theorem \ref{thm:c2}. }
Part \ref{G1} is well known and classical. We give
here a simple proof based on the analysis of foci.
Let us suppose that $X\subset\bb{P}^N$ is a 3-dimensional
variety with Gauss map whose
fibres have dimension 2. Let us consider the
family $B\subset\bb{G}(2,N)$ of the fibres of the Gauss map of $X$.
Then, by Theorem
\ref{thm:B}, there is a focal line on every general plane of $B$.
If there is a fundamental line $L$, $X$ must be a cone over a curve,
with vertex $L$. Otherwise, the focal locus is a ruled
surface $S$, and, by Theorem \ref{thm:A}, every plane of $B$ is tangent to
$S$ along a line of its ruling. Hence $S$ is a surface with
degenerate Gauss mapping, so $S$ is a cone or the tangent
developable to a curve. In the first case, $X$ is a
cone, with a point as vertex, over the tangent developable to a curve.
In the second case, $X$ is the union of osculating planes to a curve.

We will prove now  part 1.
For more details  see also \cite{tesi}.

Let $B\subset\bb{G}(1,N)$ be the
family of fibres of the Gauss map of $X$. By
\ref{thm:C}, on a general line of $B$ there are two foci (counting
multiplicity).  Then, we will consider the number of distinct focal
points on a general line belonging to $B$, the number (1 or 2) and the
dimension of the irreducible components of the strict focal locus
(see Definition \ref{d:strict}), i.e. the variety obtained as closure
of the union of focal points on non-focal lines of $B$. This is a
general procedure, which will be extended to varieties of higher
dimension in \S \ref{s:rette}.

If the focal points on a general line
of $B$ are distinct, Theorem \ref{thm:A} gives us the classification
of all possible cases, as arranged in  Table~1.
\begin{table}
\begin{tabular}{|m{4cm} | m{3cm}| m{3cm}|}
\hline\hline
foci on a general line & strict focal
locus & description\\ \hline\hline
two distinct points&$\mspace{-9mu}$\begin{tabular}{m{3cm}|m{3cm}|}
each point sweeps a
surface&$\mspace{-9mu}$\begin{tabular}{m{3cm}}{{union of lines
bitangent to a surface}}
\\\hline{{union of lines tangent to two surfaces}}\\\end{tabular}\\ \hline
  a point sweeps a surface, the other sweeps a curve &
union of lines tangent to a surface and meeting a curve\\\hline
	each point sweeps a curve&$\mspace{-9mu}$\begin{tabular}{m{3cm}}secant
variety of a curve \\\hline join  of two curves\\\end{tabular}\end{tabular}\\
\hline
\end{tabular}
\caption{Two distinct foci}
\end{table}

If on a general line of $B$ there is one focal point
with multiplicity 2, we need more information.
That can be provided considering the interpretation, given in
Section \ref{sec:foci}, of the
characteristic map of $B$
relative to a general $r\in B$ as describing
the subspace $T_rB\subset T_r\bb{G}(1,N)\cong\Hom(\hat{r},V/\hat{r})$.

Using  it,
  we will prove now that, if the strict focal locus is a surface $S$, then a
general line of $B$ represents an asymptotic direction of $S$, i.e. a
selfconjugate direction  with respect
to the second fundamental form of $S$.
Let $r$ be a general fibre of the Gauss map and $F=[v]$ be the 
(double) focus on $r$:
then by Remark \ref{fuochimult} there exist two linearly independent 
tangent vectors $\eta_1$,
$\eta_2$ in
$T_rB$ such that $\eta_1(r)=0$ in $V/\hat r$ and $\eta_2(v)\in 
Im\ (\eta_1)$. Let $\{b_1(t)\}$
be an arc of smooth curve in $B$, parametrized by an open disc 
containing the origin, with
$b_1(0)=r$ and
$b_1'(0)=\eta_1$ and let
$\{c_1(t)\}$ be a lifting of
$\{b_1(t)\}$ through $F$, i.e. any regular curve in $X$ such that 
$c_1(t)\in b_1(t)$ for all
$t$ and $c_1(0)=F$. Then $r$ is the tangent line to the curve 
$\{c_1(t)\}$ at $F$. In
particular, the curve
$C$ ($\subset S$) generated by the unique focus of $b_1(t)$ as $t$ 
varies in the disc 
is such  a lifting.  If
$\{g_1(t)\}$ is another  lifting of $\{b_1(t)\}$ but with $g_1(0)\neq 
F$, then the tangent
vector
$g_1'(0)$ is not parallel to $r$ so, together with $r$, it generates 
the tangent plane at
$g_1(0)$ to the ruled surface
$Y$, union of the lines  $b_1(t)$. Since dim $Im\ (\eta_1)=1$, this 
plane is constant along
$r$, so
it coincides with the osculating plane to the curve $C$ at $F$.
Let now $\{b_2(t)\}$ be a regular curve in $B$ such that $b_2(0)=r$ 
and $b_2'(0)=\eta_2$: if
$\{d_2(t)\}$ is a lifting of its 
through $F$, then its tangent line 
at $F$ is contained in the
osculating plane to $C$ at $F$. In particular, we can choose as 
lifting the curve $D$ of the
foci of the lines $b_2(t)$. Because of the generality assumptions, 
the tangent plane to $S$
at $F$ is generated by the tangent lines to $C$ and $D$, so it is the
osculating plane to 
$C$. We have thus proved that through a general point $F$ of $S$ there 
is a curve $C$ whose
osculating plane at $F$ coincides with the tangent plane to $S$. The 
tangent line to $C$,
which is a general line of $B$, is therefore an asymptotic tangent
line of $S$: this proves
our claim.

\begin{table}
\begin{tabular}{|m{4cm} | m{3cm}| m{3cm}|}
\hline\hline
foci on a general line & strict focal locus & description\\
\hline\hline
a point with multiplicity 2&$\mspace{-9mu}$\begin{tabular}{m{3cm}|m{3cm}|}
the focal point sweeps a surface & union of asymptotic lines\\\hline the focal
point sweeps a curve &band\\\hline
	the focal point is fixed & cone\end{tabular} \\
\hline
\end{tabular}
\caption{One double focus}
\end{table}

If the strict
focal locus is a curve $C$, then we we will show that $X$ is not just 
a union of
cones, as in the
case in which the focal point has multiplicity 1, but a 
union of planes tangent to $C$. We proceed as in the previous case: 
let $r\in B$ be a general
line and $F=[v]$ be its focus. Since $F$ is a fundamental point for 
the family $B$, there is a
curve $Z$ in the Grassmannian, passing through $r$, which represents 
the lines of $B$
through $F$. It is easy to show that every lifting of $Z$ through $F$ 
has $r$ as tangent line
at $F$, so $\eta_1$, the tangent vector to $Z$ at $r$,  is such that 
$\eta_1(v)=0$. But $F$
is a focus with multiplicity two and dim $T_rB=2$, so it follows that 
every regular curve
contained in $B$, passing through $r$ but with tangent vector  
$\eta_2$ different from
$\eta_1$, is such that $\eta_2(v)$ belongs to the image of $\eta_1$. 
The focal curve $C$ can
be interpreted as a lifting of such a curve: let $w$ be its tangent 
vector at $F$, then the
plane generated by $r$ and $w$ contains also the tangent line to any 
lifting of $Z$ at its
intersection point with $r$. Let $\phi(t)$ be a local parametrization 
of such a lifting,
with $\phi(0)=P\in r$,
then we have that $\phi'(0)$
lies in the plane generated by $w$ and the direction of $r$.
By repeated derivations,
we get that all derivatives $\phi^{(k)}(0)$ belong to this plane, 
hence the whole curve is
contained in it. Therefore every lifting of $Z$ is a plane curve, 
which proves that the lines
of $B$ passing through $F$ form a pencil, contained in the plane 
generated by $r$ and by the
tangent line to $C$ at $F$.

So  $B$ is a ruled surface.  In this
  case, $X$ is called a
3-dimensional band.
The precise definition is the
following (see \cite{AG}):
\begin{Def}\label{band}
A variety $X\subset \mathbb{P}^N$ is said to be a 
\textit{3-dimensional band} if
there exist two distinct curves $C,D\subset X$, not belonging both
to the same $\mathbb{P}^3$, and a birational equivalence $\psi: C
\longrightarrow D$, such that $X$ is the closure of the union of the
planes lying in the image of the morphism: $$\begin{array}{lccc}f:
&C_0 &\longrightarrow &\mathbb{G}  (2,N)\\&p &\longrightarrow 
&\la\mathbb{T}_pC,
\psi(p)\ra,\end{array}$$ where $C_0$ is a non-singular open
subset of $C$ contained in the domain of definition of $\psi$.\end{Def}
Table~2  describes all varieties ruled by a 2-dimensional family
of lines with a focal point of multiplicity 2 on the general line.
\qed

\vskip 12pt
By a direct calculation, it is possible to find out that every variety
obtained in Theorem \ref{thm:c1} is a variety with degenerate Gauss 
mapping. 
The interesting point
is that, whereas any curve can be obtained as the focal
curve of a $3$-dimensional variety with degenerate Gauss mapping, the
focal surfaces must verify some special conditions. For instance, it
is not a general fact for a surface that there exists a family of
dimension $2$ of bitangent lines.

\begin{thm}
Let $X\subset\bb{P}^N$ be a variety of dimension 3
with Gauss image of dimension 2, satisfying one
of the conditions (a)-(d) in Theorem \ref{thm:c2}. \\
Suppose that the strict focal
locus of the family $B$ of the fibres of the Gauss map of $X$
has an irreducible component $S$ of dimension 2.
Then $S$ is either a developable surface or a $\varPhi$ surface. Moreover  the
lines of
$B$ are tangent  to $S$ and they are either  asymptotic 
tangent lines  or they  admit a conjugate direction.
\end{thm}
\proof
If $X$ is as in (d), then the Theorem is clearly true. So we 
assume that on a general fibre
of the Gauss map there are two distinct foci. Let $F_1\in S$ be general: it is
a focus on a non-focal line $r$, which contains also a second focus $F_2$. So
there exist two tangent vectors $\eta_1, \eta_2\in T_rB$, such that, for all
regular curves 
$\{b_i(t)\}\subset B$,
$i=1,2$, with $b_i(0)=r$ and $b_i'(0)=\eta_i$, every lifting through 
$F_i$ has $r$ as tangent
line at $F_i$. As a lifting of $\{b_1(t)\}$, we can choose a curve 
$C_1\subset S$, with local 
parametrization $\{c_1(t)\}$, such that  $c_1(t)$ is a focus of the 
line $b_1(t)$ for all $t$. 
Note that $Im\ (\eta_1)$, which is 1-dimensional, is 
generated by the tangent
vector $x'(0)$ for all choice of a lifting $x(t)$ of $b_1(t)$ with 
$x(0)\neq F_1$. Hence, as
in the proof of Theorem  \ref{thm:c2}, one proves that $x'(0)$ 
belongs to the osculating
plane to the curve $C_1$.

Assume now that $X$ satisfies conditions (a) or (b). 
Then, the previous construction 
can be repeated 
for the second focus $F_2$ on $r$ 
relatively to the focal surface $S'$ to which it belongs, 
which coincides with $S$
in case (a) or is the second component of the strict focal locus of 
$X$ in case (b). This
gives a second curve $C_2\subset S'$ passing through $F_2$ and
with $\bb T _{F_1}C_1=r=\bb T _{F_2}C_2$. Let now.
  $D_2$ be the curve generated by the second focus of the lines $b_1(t)$, and
similarly $D_1$ be the curve generated by the second focus of the lines
$b_2(t)$. Note that
$C_1\neq D_1$ and
$C_2\neq D_2$. We can choose 
a local parametrization for $S$ of the form
$\psi(t,s)$, where $\psi(0,0)=F_1$, 
$\psi(t,0)$ and $\psi(0,s)$ are local parametrizations of 
respectively $C_1$ and $D_1$. 
By considering the other focus, we get a parametrization  $\phi(t,s)$ 
of the second surface $S'$ near $F_2$  
such that 
$\phi(t,0)$ and $\phi(0,s)$ are local parametrizations of 
respectively $D_2$ and $C_2$. 
By comparing the tangent vectors, we get: $\psi_t=\phi_s$,
$\phi_t\in\langle\psi_t,\psi_{tt}\rangle$,
$\psi_s\in\langle\phi_s,\phi_{ss}\rangle$, and also
$\psi_{tt}\in\langle\phi_t, \phi_s \rangle$,
$\phi_{ss}\in\langle\psi_t,\psi_s\rangle$.
So we can deduce that $\psi_{ts}\in\langle\psi_t,\psi_s\rangle$. Hence the
pair of vectors $(\psi_t, \psi_s)$ annihilates the second fundamental form of
$S$, and they represent conjugate directions. 

If we are in case (c), $F_2$ is a fundamental point for the family $B$, so
there are infinitely many lines of $B$ through $F_2$. Each of them contains
also a second focus, describing a curve $E$. In this case we can find a local
parametrization of $S$, $\psi(t,s)$, centred at $F_1$ 
and such that
$\psi(t,0)$ describes $C_1$ and
$\psi(0,s)$ describes $E$. Note that $\psi_t(0,s)$ is the direction of the
line of the ruling passing through $\psi(0,s)$, and $\psi_{ts}$ is  tangent at
$F_1$ to the cone of vertex $F_2$ on the curve $E$, therefore it is contained
in the tangent plane to this cone along $r$. But this plane is generated by
$\psi_t$ and $\psi_s$, so it coincides with the tangent plane to $S$ at the
point $F_1$. We conclude then as in the previous case. 
\qed

\vskip 12pt
We will close this section with a remark on the second fundamental 
form. It is known
(\cite{GH79}) that the second fundamental form of the varieties with degenerate
Gauss mapping has non-empty singular locus. In 
particular, this singular locus is a point in the case of varieties 
of dimension
3 with Gauss image of dimension 2. Assume that $X$ is such
a variety, which is not a hypersurface. Then there is a connection between the
properties of the second fundamental form and the configuration of
focal points on the general fibre of the Gauss map of $X$. Indeed, if $X$ is
a  variety with distinct focal points  of multiplicity 1, then the
dimension of the second osculating space is 5 and the second fundamental form
  is a pencil of conics with a point both as base  and as singular
locus. If $X$ has one  focal point of multiplicity 2 on the general line and
is  not a cone, then the dimension of the second osculating
space is also 5, but the pencil of conics of the second fundamental form has
a line as
  base locus. In the case of cones over
a surface, the dimension of the second osculating space is 6 (in general),
so the second fundamental form is a net of conics and the base locus can
only be a point, coinciding with the singular point.

\
\section{Varieties covered by lines}\label{s:rette}
Let $B\subset\mathbb{G}(1,N)$ be a family of lines in $\mathbb{P}^N$ of
dimension $n\leq N-2$. Suppose that the union of the lines belonging to $B$
is an algebraic variety $X$ of dimension $n+1$. When this
condition holds, we do not expect in general cases to find any focal point.
In particular, a general line of $B$ cannot be focal. This
allows us to consider the length of the focal locus on the general $r\in B$.
It turns out that such length has a geometric interpretation
in terms of fixed tangent spaces along $r$. Theorem \ref{thm:C} states that
  the length of the focal locus on $r\in B$ is $k$ if and only if
$X$ possesses a fixed tangent space of dimension $k+1$ along a 
general line $r$.
We will prove it now.
\vskip 12 pt
\par {\em Proof of Theorem \ref{thm:C}.}
We can suppose without loss of generality that $X$ is a hypersurface, 
i.e. $N=n+2$.
Indeed, if $X$ is not a hypersurface, we can project it to $\mathbb{P}^{n+2}$,
  and a general projection will not affect either its
tangential properties or its focal ones.
%First part
Let $r$ be a general point of $B$. Suppose that on $r$ there are exactly
$k$ focal points, counting multiplicity. They are
the points where the characteristic map relative to $r$,
$$\begin{CD}\chi(r):\ @. T_rB\otimes\mathcal{O}_r@>>>
\mathcal{N}_{r|\bb{P}^N}\\ @. \wr@| \wr@|\\ 
@.\mathcal{O}^n_r@.\mathcal{O}_r(1)^{n+1}
  \end{CD}$$
has not maximal rank. If we choose projective coordinates
$x_0,x_1$ on $r$, by means of the natural identification
given above, we can represent $\chi(r)$ by an $n\times
(n+1)$ matrix $$A=\begin{pmatrix}l_{1,1}&\dots&l_{1,n}\\
\vdots&\cdots &\vdots
\\l_{n+1,1}&\ldots&l_{n+1,n}\end{pmatrix},$$ whose
entries $l_{i,j}$ are linear forms in $x_0,x_1$. Let us
consider the minors (with sign)
of $A$ of maximal order, $$\phi_i=(-1)^{i+1}\det
\begin{pmatrix}l_{1,1}&\dots&l_{1,n}\\ \vdots&\cdots &\vdots
\\\widehat{l_{i,1}}&\dots&\widehat{l_{i,n}}\\ \vdots&\cdots &
\vdots \\l_{n+1,1}&\ldots&l_{n+1,n}\end{pmatrix},\ i=1,\dots,n+1.$$
The existence of $k$ focal points implies that $\phi_1,\dots,
\phi_{n+1}$ have a common factor $F$ of degree $k$. So we have
the relations $\phi_i=F \psi_i$, where $\psi_1,\dots,\psi_{n+1}$ are suitable
polynomials of degree $n-k$ in $x_0,x_1$. We are interested
in finding vectors tangent to $X$ at every point of $r$.
This means we seek
normal vectors of coordinates $(v_1,\dots,v_{n+1})$
belonging to the image of $\chi(r)$ in every point of $r$. This can 
be expressed by
the condition
$$\det\begin{pmatrix}v_1&l_{1,1}&\dots&l_{1,n}\\
\vdots&\vdots&\cdots &\vdots \\v_{n+1}&l_{n+1,1}&
\ldots&l_{n+1,n}\end{pmatrix}=0,$$
or, equivalently, \begin{equation*}\tag{*}v_1\psi_1+ 
\dots+v_{n+1}\psi_{n+1}=0.\end{equation*} As there are
$n-k+1$ monomials of degree $n-k$ in $x_0,x_1$, equation
(*) is equivalent to a system of $n-k+1$ homogeneous linear
equations
in the indeterminates $v_1,\dots,v_{n+1}$. So there are at
least $k$ linearly independent solutions. Denote by $V'$
a linear subspace of dimension $k$ of the space of solutions. If we identify
$V/\hat{r}$ with a subspace $W$ complementary of $\hat{r}$,
the vectors of $V'\subset V/\hat{r}$ are tangent to $X$ at  every point of
$r$. Then $\mathbb{P}(V')$ is a tangent subspace of
dimension $k+1$ contained in the tangent space to $X$ at every point 
of $r$. That
proves implication (i) $\Rightarrow$ (ii).
%Second part
Let $r$ be a general point of $B$. Now we will prove that
if there is a constant tangent space of dimension $k+1$ along
$r$ then there are $k$ focal points on $r$ (counting multiplicity).
As in the previous part, we will denote by $A=(l_{i,j})$ the matrix
representing $\chi(r)$. What we want to show is that the minors
$\phi_1,\dots,\phi_{n+1}$ have a common factor of degree $k$.
We know that condition $$\det\begin{pmatrix}v_1&l_{1,1}&\dots&l_{1,n}\\
\vdots&\vdots& \cdots &\vdots 
\\v_{n+1}&l_{n+1,1}&\ldots&l_{n+1,n}\end{pmatrix}=0$$
is satisfied for every $v=(v_1,\dots,v_{n+1})$
belonging to a normal subspace of dimension $k$.
We can assume without loss of generality that this normal subspace is
$$V'=\la(\overbrace{0,\dots,0}^{n-k+1},1,\overbrace{0,
\dots,0}^{k-1}),( \overbrace{0,\dots,0}^{n-k+2},1,\overbrace{0,\dots,0}^{k-2}),
\dots,(\overbrace{0 ,\dots,0}^n,1)\ra.$$ This is the
same as supposing that the last $k$ minors of order $n$ of $A$ are 0, i.e.
$$\phi_{n-k+2}=\dots=\phi_{n+1}=0.$$
In the following, we will denote by $A^{j_1,\dots,j_h}_{i_1,\dots,i_h}$
the determinant of the square submatrix of the $j_1,\dots,j_h$-th rows
and the $i_1,\dots,i_h$-th columns of $A$.
Let us consider the remaining forms $\phi_{1},\dots,\phi_{n-k+1}$.
Being minors of the matrix $A$, they satisfy the following 
homogeneous system of degree 1
$$\sist
l_{1,1} \phi_1 + l_{2,1} \phi_2 + \dots + l_{n-k+1,1}
\phi_{n-k+1} = 0\\ \hdotsfor{1}\\
l_{1,n} \phi_1 + l_{2,n} \phi_2 + \dots + l_{n-k+1,n} \phi_{n-k+1}
= 0.\sistt$$ Fix two equations of the system above by
choosing two
  indices $1\leq i_1<i_2\leq n$.
We can multiply the first equation by $l_{n-k+1,i_2}$,
  the second one by $l_{n-k+1,i_1}$ and subtract: we
get a homogeneous relation 
among $\phi_1,\dots,\phi_{n-k}$ with coefficients of degree 2.
Considering every possible choice of $i_1,i_2$,
we obtain the homogeneous system
$$\sist 
A^{1,n-k}_{i_1,i_2}\phi_1+\dots+A^{n-k,n-k+1}_{i_1,i_2}\phi_{n-k}=0\\ 
1\leq i_1<i_2\leq n.\sistt
$$
In an analogous way we can find homogeneous relations with coefficients
of every degree between 2 and $n-k$, involving less and less minors.
For the highest degree we have a system of $\binom{n}{k}$
equations in 2 minors. For $\phi_1,\phi_2$, for instance, we have
$$\sist
A^{1,3,4,\dots,n-k+1}_{i_1,\dots,i_{n-k}}\phi_1+ 
A^{2,3,\dots,n-k+1}_{i_1,\dots,i
_{n-k}}\phi_2=0\\
0\leq i_1<i_2<\dots<i_{n-k}\leq n.\sistt
$$
As the relations belonging to this system cannot all be trivial,
we get that $\phi_1$ and $\phi_2$ must have a common factor of
degree $\geq k$. Moreover, it is possible to prove that
$\phi_1,\dots,\phi_{n-k}$ have a common factor of degree $k$ .
In fact, consider the system of relations (of degree $n-k-1$)
among 3 minors, say $\phi_1,\phi_2,\phi_3$:

$$\sist A^{1,4,\dots,n-k+1}_{i_1,\dots,i_{n-k-1}}\phi_1+A^{2,4,
\dots,n-k+1}_{i_1, \dots,i_{n-k-1}}\phi_2+A^{3,4,\dots,n-k+1}_{i_1,
\dots,i_{n-k-1}}\phi_3=0\\ 0\leq i_1<i_2<\dots<i_{n-k-1}\leq n.\sistt
$$
Denote by $F$ a common factor of degree $k$ of $\phi_1,\phi_2$.
Suppose $F\nmid \phi_3$: then there is a factor $G$ of $F$ such
that $G$ divides $A^{3,4,\dots,n-k}_{i_1,\dots,i_{n-k-2}}$ for
any choice of $i_1,\dots,i_{n-k-2}$. This means that $G$ divides
both $A^{1,3,4,\dots,n-k}_{i_1,\dots,i_{n-k-1}}$ and
$A^{2,3,\dots,n-k}_{i_1,\dots,i_{n-k-1}}$. Then $\phi_1$ and $\phi_2$
have a common factor of degree $\geq k+\deg G$, and we can check
whether this new polynomial of higher degree and $\phi_3$
have a common factor of degree $k$ or not. If the answer is negative,
we can iterate the construction until we find the factor we look for,
after deleting all common factors of
$A^{1,3,4,\dots,n-k}_{i_1,\dots,i_{n-k-1}}$ and $A^{2,3,
\dots,n-k}_{i_1,\dots,i_{n-k-1}}$.
\qed
\vskip 12pt
\begin{rem}
If there are more than $n$ focal points on a line $r\in B$,
then $r$ is a focal line.
\end{rem}
Theorem \ref{thm:C} allows us to give a rough description
of the focal locus of a variety $X$ ruled by an $n$-dimensional family of
lines, once we
know the dimension of the constant tangent space along a general
line.
\begin{Def}\label{d:strict}
Let $B\subset\bb{G}(1,N)$ be a subvariety of the Grassmannian,
such that its general element is not focal. Let $k$ be the
degree of the focal locus on a general element $r\in B$.
Let us denote by $U\subset B$ the open set of the lines on which the 
focal locus
is a proper subscheme of length $k$. Then the closure in
$\bb{P}^N$ of the union of the focal loci on the lines of $U$ is 
called the {\em
strict focal locus} of $B$.\end{Def}
\begin{rem} For $B\subset\bb{G}(1,N)$, $\dim B=n$, if $n=k$,
then $U$ is the open set of non-focal lines, and the strict focal
locus is simply the closure in $\bb{P}^N$ of the union of
focal points on the non-focal lines of $B$.\end{rem}
The study of the strict focal locus enables us to formulate a
pattern of classification of varieties of dimension $n+1$ ruled
by an $n$-dimensional family of lines. First of all, any such variety is
characterized by the number and the multiplicity of the distinct
focal points on a general line. Then we can study the strict focal locus
of $X$, and, in particular, the number of components  and their 
dimensions.\\
The maximal possible dimension for a component of the strict
focal locus is $n$. If there is a component of dimension $<n$,
then through every point of it there pass 
infinitely many lines
 of $B$. So
this component must be contained in the fundamental locus of the
family $B$. If there are components of dimension $n$ of the focal locus,
then every line of $B$ is tangent to them. This is a particular
  case of a property of the focal locus that holds for varieties ruled by
linear  subspaces of dimension $\geq 1$ too. So we will prove
it in the general case.
\begin{thm}\label{thm:A}
Let $B\subset\bb{G}(k,N)$ be a family of $k$-spaces in $\bb{P}^N$
of dimension $n\leq N-k$. Suppose that the union of the $k$-planes 
belonging to $B$ is a
variety $X$ of dimension $k+n$, and that the focal
locus has codimension 1 in $X$.\\ Then every general subspace 
$\Lambda$ belonging to
$B$ is tangent to $F$ at all the focal points on
$\Lambda$ that are not fundamental points.
\end{thm}
\proof
Let us consider $\mathcal{I}$, the  desingularization of the
incidence correspondence of $B$, and the natural projections $f,g$
$$\begin{CD}\mathcal{I}@>f>>\mathbb{P}^N\\
@VgVV\\
B.\end{CD}$$
Let $p$ be a general point of $g^{-1}(\Lambda)$, belonging
to the focal subvariety $V(\chi)\subset\mathcal{I}$. 
By definition, the differential 
of $f$ in $p$,
$$d_pf: T_p\mathcal{I} \longrightarrow T_{f(p)} \bb{P}^N$$
has a non-trivial kernel.\\
Let us consider the restriction of $f$ to the focal subvariety,
$$\begin{array}{lccc}f: & \mathcal{I}& \xrightarrow{\ \ \ \ \ \ \ 
}&\,\;\mathbb{P}^N\\[-3pt]
&\scriptstyle{\cup} &&\scriptstyle{\cup} \\
&V (\chi)& \xrightarrow{\ \ \ \ \ \ \ }& F.\end{array}$$
As we
supposed $char\ \bb{K}=0$, the algebraic geometry analogue of
Sard's Theorem holds (\cite [p. 176]{H92}). So 
for a general $p\in 
V( \chi)$ the homomorphism 
$d_pf|_{V (\chi)}: T_pV(\chi) \longrightarrow T_{f (p)}F$ is
surjective.
	Since not all focal points are fundamental points and $f$ has
finite-dimensional fibres, $\dim V(\chi)=\dim F$. Thus $d_pf|_{V(\chi)}$
is an isomorphism.\\ Now consider again the differential of $f$
in a general point $p$ of $V(\chi)\cap g^{-1}(\Lambda)$, $$d_pf:
T_p\mathcal{I} \rightarrow T_{f (p)} \mathbb{P}^N.$$
We know that $d_pf$ has a non-trivial kernel. We already know as
  well that its image contains $T_{f (p)}F$. Since $V (\chi)$ is a
codimension 1 subvariety of $\mathcal{I}$, $T_pV (\chi)$ is a linear
subspace of codimension 1 in $T_p{\mathcal{I}}$. Hence
$d_pf (T_p\mathcal{I}) = T_pF$.  As $g^{-1} (
\Lambda) \subset \mathbb{T}_p\mathcal{I}$, we
have $\Lambda = d_pf (g^{-1} (\Lambda) \subset \mathbb{T}_{f (p)}F$.
\qed
\vskip 12pt
Coming back to varieties with constant tangent space along lines,
assume that on a general line $r\in B$ there are focal
points which are not fundamental points. Then the strict focal
locus has a
component $Y$ of pure codimension 1 in $X$, and every line in $U$ is
tangent to $Y$ at its focal, non-fundamental points.
In Theorem \ref{thm:C} we can consider the two extremal cases:
namely, $k=0$ and $k=n$. In the first case, the theorem  implies that
  a variety ruled by lines has no focal point on a general line if and only
if the only fixed tangent space along a general line is the line
itself. In the second case, we obtain a characterization of the varieties
whose degenerate Gauss map has 1-dimensional fibres, i.e.
varieties of dimension $n+1$ with tangent space constant along lines.
\begin{cor}\label{cor:C}
Let $B\subset\mathbb{G}(1,N)$ be a family of lines in $\mathbb{P}^N$
  of dimension $n,\ n\leq N-2$. Suppose that the union of the lines
belonging to
$B$ is an algebraic variety $X$ of dimension $n+1$.
Then 
the following are equivalent:\\
(i) the focal locus on the general element $r\in B$ consists of $n$ points
(counting multiplicity);\\
(ii) the tangent space to $X$ is constant along the lines of $B$. \end{cor}
\begin{rem}
If $X$ is not ruled by linear subspaces of dimension $\geq2$,
then condition (i) implies that $B$ is the family of the fibres of 
the Gauss map.
If $X$ possesses a higher dimensional ruling, then the
fibres of the Gauss map may have dimension greater than 1.\end{rem}

\section{Varieties ruled by linear subspaces}\label{gen}
In this section we try to find out whether the results established in the
previous section may be extended to 
varieties ruled by linear subspaces of dimension $>1$.
In particular we expect that in the case of a 
family of linear subspaces of dimension $k$, 
the existence of constant tangent
spaces gives a focal hypersurface on the general $k$-space.
This is true for varieties with degenerate Gauss mapping, 
for which
Theorem \ref{thm:B} yields a straightforward generalization 
of Corollary \ref{cor:C}.
\vskip 12pt
\par{\em Proof of Theorem \ref{thm:B}.}
%%%first part
Let $X\subset\mathbb{P}^N$ be a projective variety of dimension $n+k$,
with Gauss map whose fibres have dimension $k$. We want to prove
that condition (ii) holds for the family $B\subset\mathbb{G}(k,N)$ of
fibres of the Gauss map of $X$. Let $\Lambda$ be a general element
  of $B$. Let us consider the characteristic map of $B$ relative to
$\Lambda$,
$$\begin{CD}
\chi(\Lambda):\ @. T_\Lambda B \otimes \mathcal{O}_\Lambda @>>>
\mathcal{N}_{\Lambda|\mathbb{P}^N}\\
@. \wr@| \wr@|\\
@. \mathcal{O}^{n}_\Lambda @>>> \mathcal{O}^{N-k}_\Lambda (1).\end{CD}$$
$\chi(\Lambda)$ is represented by an $n\times(N-k)$ matrix,
whose entries are linear forms
on $\Lambda$. The columns of this matrix evaluated in a point
$p\in\Lambda$ can be regarded as vectors $L_1 (p),\dots,L_{n}(p)$
in $V/\hat{\Lambda}$. Let us denote by $\Pi$ the fixed tangent space
to $X$ along $\Lambda$. Then the image of $\chi(\Lambda)$ in any
point $p\in\Lambda$ is contained in $\hat{\Pi}/\hat{\Lambda}$, which
is a fixed subspace of $V/\hat{\Lambda}$ of dimension $n$. If we
consider the coordinates of the normal vectors $L_1 (p),\dots,L_{n}(p)$
in $\hat{\Pi}/\hat{\Lambda}$, we get a matrix $(m_{ij})_{i,j=1,\dots,n}$.
Then the condition defining the focal locus on $\Lambda$ is 
$$\det(m_{ij}) = 0,$$
which in general cases gives a hypersurface on $\Lambda$
of degree $n$, even if it is possible in special cases that 
all $\Lambda$ is focal.
%%%second part
Suppose now that $B\subset\mathbb{G}(n,k)$ satisfies condition
(ii). If we fix a general $\Lambda\in B$, the focal variety on $\Lambda$
is a hypersurface of degree $n$. Then on the general line $r\subset\Lambda$
there are $n$ (not necessarily distinct) focal points, which
are the points where the morphism
$$\lambda:\ T_\Lambda B\otimes\mathcal{O}_r 
\longrightarrow(\mathcal{N}_{\Lambda|\mathbb{P}^N})|_r,$$
given by the restriction of the characteristic
map, has not  maximal rank. We can adapt to $\lambda$ the procedure 
applied in the proof of
the implication (ii) $\Rightarrow$ (i) of Theorem 
\ref{thm:C}. In this way, we find that there is a fixed subspace 
$W(r)$ of dimension
$n$ contained in the image of the characteristic map
in any point of $r$. Now we choose a general point $p$ in $\Lambda$; 
particularly, $p$
is non-focal and smooth. We restrict to an affine open
set $U_0\subset\mathbb{P}^N$ and consider a system of
affine coordinates on $U_0$
such that $p$ is the point $(0,\dots,0)$. $\Lambda_0=\Lambda\cap U_0$
is a linear space of dimension $k$. We can fix $k$ lines
$r_1,\dots,r_k$ through $p$\ spanning $\Lambda_0$, such that
  for all $j$ on $r_j$ there are $n$ focal points (considered with
multiplicity). On any $r_j$ there is a fixed tangent subspace,
  spanned by $r_j$ and $W(r_j)$. So the tangent space to $X$ in $p$ does
contain all lines $r_1,\dots,r_k$ (spanning $\Lambda_0$) and
all linear subspaces $W(r_1),\dots,W(r_k)$, which implies that all the spaces
$W(r_j)$ must coincide for dimensional reasons. In this way we
have found a fixed linear space $W$ of dimension $n$, such that in any
smooth point of $\Lambda_0$ the tangent space to $X$ is spanned by 
$\Lambda_0$ and $W$. \qed
\vskip 12pt
In the general case of varieties ruled by lines, it was possible
to find non-focal lines on which there were more than
the general number of focal points. Under the hypotheses of  
Theorem  \ref{thm:B}, the
open set of subspaces on which the focal locus has
  degree $k$ coincides with the set of non-focal subspaces.
So Theorem \ref{thm:B}
allows us to describe possible characterizations
  of the strict focal locus for varieties with degenerate Gauss mapping.
In this case, the strict focal locus is defined
  as the closure in $\bb{P}^N$ of the union of the focal points on 
non-focal subspaces.
For varieties with degenerate Gauss mapping, the focal locus
is contained in the singular locus of the variety.
The converse is not true in general.

\begin{thm}
Let $X$ be a variety with degenerate Gauss mapping, and denote by
$B$ the family of fibres of the Gauss map of $X$. Then,
the focal points of $B$ are singular points of $X$.
\end{thm}
\proof Let us recall that the focal points are the
ramification points of the projection $f: \mathcal{I}\longrightarrow X$
  from the desingularization of the incidence correspondence of $B$
to $X$. As the degree of $f$ is 1, either
  the focal points are points where $f$ is not
finite, or they are necessarily non-normal points of $X$. In the
former case, they are fundamental points of $B$;
in the latter, they are a fortiori singular points of $X$. Since
through a fundamental
point there pass at least two different fibres of the Gauss map,
  also the fundamental points are always singular.
\qed

\begin{rem}
This theorem can be extended to every variety $X$ ruled by a
family $B$, such that the projection from the desingularized
incidence corrispondence to $X$ has degree 1. In this case all
  focal points not belonging to the fundamental locus are non-normal points,
but nothing can be said about fundamental points. \end{rem}
Theorem \ref{thm:B} could suggest that
also a more general equivalence holds true, i.e. that, given
a variety $X$ of dimension $n+k$ covered by a family $B$ of $k$-spaces with
$\dim B=n$, $X$ possesses a constant tangent space of dimension
$k+h$ along a general $\Lambda\in B$ if and only if the focal locus of $B$
on a general $\Lambda\in B$ is a hypersurface of degree $h$.
Unfortunately, there are counterexamples of this equivalence even for the
first possible non-trivial case, that is for varieties
ruled by a family of planes with focal lines. Observe that this case is the
simplest possible not covered by Theorem \ref{thm:C} or Theorem 
\ref{thm:B} either.

\begin{ese}We will give two examples of varieties of dimension
4 ruled by a $2$-dimensional family $B$ of  planes, with a focal line
on the general $\Lambda\in B$. We will see that the tangential
properties along the planes of the ruling are not the same in the two cases.

Let us  consider a variety $Y$ of dimension 3 ruled by lines,
with a fixed tangent plane along the general line
of the ruling, but no higher dimensional constant tangent space. Then 
the family of
tangent planes has a focal line on the general element,
and this line is precisely the line of the ruling of $Y$. In this case it is
possible to prove that the union of the family of
  tangent planes is  a variety $X$ of dimension 4 with a fixed tangent
$\bb{P}^3$ along every plane. So, for the variety $X$ the relationship between
  the dimension of the fixed tangent space along
the planes of the ruling and the degree of the focal locus holds.

Now, let $Z$ be a variety of dimension 3
  ruled by lines, with constant tangent
space along the lines of the ruling. Denote by $B$ the  2-dimensional 
family of such lines, i.e. (in
general) the family of the fibres of the
  Gauss map. Then we can choose a family $C\subset\bb{G}(2,N)$ of 
planes such that, for every line
$r$ in $B$, there is a plane in $C$ containing
$r$ and lying in the constant tangent space to $Z$ along $r$. On a 
general plane $\Pi$ in
$C$ the line $r$ of $B$ such that $r\subset\Pi$ is a
  focal line. Assume that the union of the planes of $C$ is a variety 
$X$ of dimension
4. It is possible to prove that along a general line in
  $\Pi$ there is a constant $\bb{P}^3$ tangent to $X$, but that this $\bb{P}^3$
depends on the chosen line, so that there is no constant
$\bb{P}^3$ tangent to $X$ along $\Pi$. This example shows 
therefore  that the relationship
previously proposed is not always valid.
\end{ese}
Concluding, all we know in general cases is that if a
  variety $X$ of dimension $n+k$, ruled by an
$n$-dimensional family of $k$-spaces,  possesses
a fixed space of dimension $k+h$ tangent
along a general $k$-space, then the focal locus
on the general $k$-space of the ruling must contain a hypersurface of 
degree $\geq h$. If
we know the degree of the focal locus, we
only know the maximal  dimension of a space tangent to $X$ along the 
general line lying in a
space of the ruling, which can vary with the choice of the line.

\end{document}